\documentstyle{amsppt}
\topmatter
\title NON cOMPLETE AFFINE CONNECTIONS ON FILIFORM LIE ALGEBRAS
\endtitle
\author Elisabeth REMM and Michel GOZE 
\endauthor
\address
Universit\'e de Haute Alsace.
Facult\'e des Sciences et Techniques.
4, rue des FR\èeres Lumi\ère F. 68093 MULHOUSE Cedex 
\endaddress
\email E.Remm\@univ-mulhouse.fr  \endemail
\email M.Goze\@univ-mulhouse.fr  \endemail

\date
\enddate
\abstract
We give a example of non nilpotent faithful representation of a filiform Lie
algebra. This gives one counter-example of the conjecture saying that every
affine connection on a filiform Lie group is complete.
\endabstract
\endtopmatter
\document
\head 1. Affine connection on a nilpotent Lie algebra \endhead

\subhead 1.1. Affine connection on nilpotent Lie algebras \endsubhead

\proclaim{Definition 1}
Let $\frak{g}$ be a $n$-dimensional Lie algebra over $\Bbb{R}$. It is called
affine if there is a bilinear mapping 
$$
\nabla :\frak{g}\times {\frak{g}}\rightarrow \frak{g}
$$
satisfying 
$$
\left\{ 
\aligned
& 1)\text{\quad }\nabla \left( X,Y\right) -\nabla \left( Y,X\right) =\left[
X,Y\right]  \\ 
& 2)\text{\quad }\nabla \left( X,\nabla \left( Y,Z\right) \right) -\nabla
\left( Y,\nabla \left( X,Z\right) \right) =\nabla \left( \left[ X,Y\right]
,Z\right) 
\endaligned
\right. 
$$
for all $X,Y,Z\in \frak{g.}$
\endproclaim

If $\frak{g}$ is affine, then the corresponding connected Lie group $G$ is
an affine manifold such that every left translation is an affine isomorphism
of $G$. In this case, the operator $\nabla $ is nothing that the connection
operator of the affine connection on $G$.

Let $\frak{g}$ be an affine Lie algebra. Then the mapping 
$$
f:  \frak{g}  \rightarrow End(\frak{g})
$$
defined by 
$$
f(X)(Y)=\nabla (X,Y) 
$$
is a linear representation (non faithful) of $\frak{g}$ satisfying 
$$
f(X)(Y)-f(Y)(X)=\left[ X,Y\right] \qquad (*) 
$$

\proclaim{Remark}
The adjoint representation $\widetilde{f}$ of $\frak{g}$ satisfies 
$$
\widetilde{f}(X)(Y)-\widetilde{f}(Y)(X)=2\left[ X,Y\right] 
$$
and cannot correspond to an affine connection.
\endproclaim

\subhead 1.2. Classical examples of affine connection \endsubhead

i) Let $\frak{g}$ be the $n$-dimensional abelian Le algebra. Then the
representation 
$$
f:  \frak{g} \rightarrow  End(\frak{g})
$$
given by 
$$  
 X  \mapsto  f(X)=0
$$
defines an affine connection.

ii) Let $\frak{g}$ be an $2p$-dimensional Lie algebra endowed with a
symplectic form : 
$$
\theta \in \Lambda ^{2}\frak{g}^{*}\text{ such that }d\theta =0 
$$
with 
$$d\theta (X,Y,Z)=\theta (X,\left[ Y,Z\right] )+\theta (Y,\left[
Z,X\right] )+\theta (Z,\left[ X,Y\right] ).$$
 For every $X\in \frak{g}$ we
can define an unique endomorphism $f_{X}$ by 
$$
\theta (adX(Y),Z)=-\theta (Y,f_{X}(Z)). 
$$
Then $\nabla \left( X,Y\right) =$ $f_{X}(Y)$ is an affine connection.

iii) Following the work of Benoist $[1] ,$ we know that
exists nilpotent Lie algebra without affine connection.

\subhead 1.3. Faithful representations associated to an affine connection \endsubhead

Let $\nabla $ be an affine connection on $n$-dimensional Lie algebra $\frak{g%
}$. Let us consider the $(n+1)$-dimensional linear representation given by 
$$
\rho : \frak{g}  \rightarrow   \frak{g}\oplus \Bbb{R} $$
given by 
$$
\rho (X):(Y,t)\mapsto \left( f_{X}(Y)+tX,0\right) 
$$
It is easy to verify that $\rho $ is a faithful representation of dimension $%
n+1$ if and only if $f_{X}(Y)=\nabla (X,Y)$ is an affine connection.

\proclaim{Definition 2}
We say that the representation $\rho $ is nilpotent if the endomorphism $%
\rho (X)$ is nilpotent for every $X$ in $\frak{g.}$
\endproclaim

\proclaim{Proposition 1}
Suppose that $\frak{g}$ is a complex nilpotent Lie algebra and let $\rho $
be a faithful representation of $\frak{g.}$ Then there exists a faithful
nilpotent representation.
\endproclaim

Proof: Let us consider the $\frak{g}$-module $M$ associated to $\rho .$\
Then, as $\frak{g}$ is nilpotent, $M$ can be decomposed as 
$$
M={\oplus }_{i=1}^{k}M_{\lambda _{i}}
$$
where $M_{\lambda _{i}}$ is a $\frak{g}$-submodule, and the $\lambda _{i}$
are linear forms on $\frak{g}$. For all $X\in \frak{g}$, the restriction of $\rho \left( X\right)$
to $M_{i}$ as
the following form 
$$
\left( 
\matrix
\lambda _{i}(X) & * & \cdots  & * \\ 
0 & \ddots  & \ddots  & \vdots  \\ 
\vdots  & \ddots  & \ddots  & * \\ 
0 & \cdots  & 0 & \lambda _{i}(X)
\endmatrix
\right) 
$$
Let $\Bbb{C}_{\lambda _{i}}$ be the one dimensional $\frak{g}$-module
defined by 
$$
\mu : X \in \frak{g} \rightarrow   \mu(X) \in End\Bbb{C} $$
with
$$ 
\mu(X)(a) = \lambda _{i}(X)a
$$
The tensor product $M_{\lambda _{i}}\otimes \Bbb{C}_{-\lambda _{i}}$ is the $%
\frak{g}$-module associated to 
$$
X\cdot \left( Y\otimes a\right) =
\rho \left( X\right) (Y)\otimes a-Y\otimes \lambda _{i}(X)a
$$
Then $\widetilde{M}=\oplus \left( M_{\lambda _{i}}\otimes C_{-\lambda
_{i}}\right) $ is a nilpotent $\frak{g}$-module. Let us prove that $%
\widetilde{M}$ is faithful. Recall that a representation $\rho $ of $\frak{g}
$ is faithful if and only if $\rho (Z)\neq 0$ for every $Z\neq 0\in Z(\frak{g%
}).$ Consider $X\neq 0\in Z(\frak{g)}.$ If$\quad \widetilde{\rho }(X)=0$,
the endomorphism $\rho (X)$ is diagonal. Suppose that $\frak{g} \neq Z(\frak{g)}$ 
and let $\Cal{C}^{k-1}(\frak{g})=Z(\frak{g})$ where $k$ is
the index of nilpotence of $\frak{g.}$ Then 
$$
\exists \left( Y,Z\right) \in \left( \Cal{C}^{k-2}(\frak{g}),\frak{g}%
\right) \diagup \left[ Y,Z\right] =X\frak{\quad }
$$
The endomorphism $\rho (Y)\rho (Z)-\rho (Z)\rho (Y)$ is nilpotent and the eigenvalues of 
$\rho (X)$ are $0$. Thus $\rho (X)=0$ and $\rho $ is not faithful. We can conclude that  $%
\widetilde{\rho }(X)\neq 0$ and $\widetilde{\rho }$ is a faithful
representation.

\head 2. Affine connection on filiform Lie algebra \endhead

\subhead 2.1. Definition \endsubhead

\proclaim{Definition 3} 
A $n$-dimensional nilpotent Lie algebra $\frak{g}$ is called filiform if the
smallest $k$ such that $\Cal{C}^{k}\frak{g}=\left\{ 0\right\} $ is equal
to $n-1.$
\endproclaim

In this case the descending sequence is 
$$
\frak{g}\supset \Cal{C}^{1}\frak{g}\supset \cdots \supset \Cal{C}%
^{n-2}\frak{g}\supset \left\{ 0\right\} =\Cal{C}^{n-1}\frak{g} 
$$
and we have 
$$
\left\{ 
\aligned
& \dim \Cal{C}^{1}\frak{g}=n-2,\text{ } \\ 
& \dim \Cal{C}^{i}\frak{g}=n-i-1,\text{ for }i=1,...,n-1.
\endaligned
\right. 
$$

\proclaim{Example}
The $n$-dimensional nilpotent Lie algebra $L_{n}$ defined by 
$$
\left[ X_{1},X_{i}\right] =X_{i+1}\text{ for }i\in \left\{ 2,...,n-1\right\} 
$$
is filiform.
\endproclaim

We can note that any filiform Lie algebra is a linear deformation of $L_{n}$ 
$[6] $.

\subhead 2.2. On the nilpotent affine connection \endsubhead

Let $\frak{g}$ be a filiform affine Lie algebra of dimension $n$, and $\rho $
be the $(n+1)$-dimensional associated faithful representation. Let $M=\frak{g%
}\oplus \Bbb{C}$ be the corresponding complex $\frak{g}$-module. As $\frak{g}
$ is filiform, its decomposition has the following form 
$$
\aligned
& 1)\text{ }M=M_{0}\quad \text{and }M\text{ is irreducible,} \\ 
& 2)\text{ }M=M_{0}\oplus M_{\lambda }.
\endaligned
$$

For a general faithful representation, let us call characteristic the
ordered sequence of the dimensions of the irreducible submodules. In the
filiform case we have $c(\rho )=(n+1)$ or $(n,1)$.\ In fact, the filiformity of $\frak{g}$ implies that exists an irreducible
submodule of dimension greater than $n-1.$ More generally, if the
characteristic sequence of a nilpotent Lie algebra is equal to $%
(c_{1},..,c_{p},1)$ (see [6]) then for every faithful representation $\rho $
we have $c(\rho )=(d_{1},..,d_{q})$ with $d_{1}\geq c_{1}.$

\proclaim{Theorem 1}
Let $\frak{g}$ be the filiform Lie algebra $L_{n}.$ There are faithful $%
\frak{g}$-modules which are not nilpotent.
\endproclaim
Proof: Consider the following representation given by the matrices $\rho
(X_{i})$ where $\{X_{1},..,X_{n}\}$ is a basis of $\frak{g}$%
$$
\rho (X_{1})=\left( 
\matrix
a & a & 0 & \cdots & \cdots &  &  & 0 & 1 \\ 
a & a & 0 &  &  &  &  & \vdots & 0 \\ 
0 & 0 & 0 &  &  &  &  & 0 & 0 \\ 
\vdots & \ddots & \frac{1}{2} & \ddots &  &  &  & \vdots & 0 \\ 
\vdots &  & \ddots & \ddots & \ddots &  &  & \vdots & 0 \\ 
\vdots &  &  & \ddots & \frac{i-3}{i-2} & \ddots &  & \vdots & 0 \\ 
0 & 0 &  &  & \ddots & \ddots & \ddots & \vdots & 0 \\ 
\alpha & \beta & 0 & \cdots & \cdots & 0 & \frac{n-3}{n-2} & 0 & 0 \\ 
0 & 0 & 0 & 0 & 0 & 0 & 0 & 0 & 0
\endmatrix
\right) 
$$
$$
\rho (X_{2})=\left( 
\matrix
a & a & 0 & \cdots & \cdots &  & \cdots & 0 & 0 \\ 
a & a & 0 &  &  &  &  & \vdots & 1 \\ 
-1 & 1 & 0 &  &  &  &  & 0 & 0 \\ 
0 & 0 & \frac{1}{2} & \ddots &  &  &  & \vdots & 0 \\ 
\vdots &  & \ddots & \ddots & \ddots &  &  & \vdots & 0 \\ 
\vdots &  &  & \ddots & \frac{1}{i-2} & \ddots &  & \vdots & 0 \\ 
0 & 0 &  &  & \ddots & \ddots & \ddots & \vdots & 0 \\ 
\beta & \alpha & 0 & \cdots & \cdots & \cdots & \frac{1}{n-2} & 0 & 0 \\ 
0 & 0 & 0 & 0 & 0 & 0 & 0 & 0 & 0
\endmatrix
\right) 
$$
and for $3 \leq j \leq n-1$ the endomorphismes $\rho (X_{j})$ satisfy : 
$$
\left\{
\aligned
&  \rho (X_{j})(e_{1})=-\frac{1}{j-1}e_{j+1} \\ 
&  \rho (X_{j})(e_{2})=\frac{1}{j-1}e_{j+1} \\ 
&   \rho (X_{j})(e_{3})=\frac{1}{j(j-1)}e_{j+2} \\ 
&  \cdots \\ 
 &  \rho (X_{j})(e_{i-j+1})=\frac{(j-2)!(i-j-1)!}{(i-2)!}e_{i},\quad i=j-2,...,n
\\ 
&  \rho (X_{j})(e_{i-j+1})=0,\quad i=n+1,..,n+j-1 \\ 
&  \rho (X_{j})(e_{n+1})=e_{j}
\endaligned
\right. 
$$
and for $j =n$  
$$
\left\{
\aligned
&  \rho (X_{n})(e_{i})=0 \quad i=1,...,n \\ 
&  \rho (X_{n})(e_{n+1})=e_n \\ 
\endaligned
\right. 
$$
where $\{e_{1},...,e_{n},e_{n+1}\}$ is the basis given by $e_{i}=(X_{i},0)$
and $e_{n+1}=(0,1).$ We easily verify that these matrices describe a non
nilpotent faithful representation.

\subhead 2.3.Study of an associated connection \endsubhead

The previous representation is associated to an affine connection on the
filiform Lie algebra $L_{n}$ given by 
$$
\nabla _{X_{i}}=\rho (X_{i})\mid _{\frak{g}}
$$
where $\frak{g}$ designates the $n-$dimensional first factor of the $(n+1)-$
dimensional faithful module. This connection is complete if and only if the
endomorphismes $R_{X}\in End(\frak{g})$ define by 
$$
R_{X}(Y)= \nabla _{Y}(X)
$$
are nilpotent for all $X\in \frak{g}$ ([5]). But the matrix of $R_{X_{1}}$
has the form : 
$$
\left( 
\matrix
a & a & 0 & \cdots  & 0 & \cdots  & 0 & 0 \\ 
a & a &  &  & \vdots  &  & \vdots  & 0 \\ 
0 & -1 &  &  & \vdots  &  & \vdots  & 0 \\ 
0 & 0 & -\frac{1}{2} & \cdots  & 0 & \cdots  & 0 & 0 \\ 
\vdots  & \vdots  & 0 & \ddots  &  & \cdots  & \vdots  & 0 \\ 
0 & 0 & \vdots  & \ddots  & -\frac{1}{j-1} &  & \vdots  & 0 \\ 
\alpha  & \beta  & \vdots  & \cdots  & \ddots  & \ddots  & 0 & 0 \\ 
0 & 0 & 0 & 0 & 0 & 0 & -\frac{1}{n-2} & 0
\endmatrix
\right) 
$$
Its trace is $2a$ and for $a\neq 0$ it is not nilpotent. We have proved :

\proclaim{Proposition 2}
There exist affine connexions on the filiform Lie algebra $L_{n}$ which are
non complete.
\endproclaim

\proclaim{Remark}
The most simple example is on $dim3$ and concerns the Heisenberg algebra. We
find a nonnilpotent faithful representation associated to the noncomplete
affine connection given by :
$$
\nabla _{X_{1}}=\left( 
\matrix
a & a & 0 \\ 
a & a & 0 \\ 
\alpha  & \beta  & 0
\endmatrix
\right) ,\quad \nabla _{X_{2}}=\left( 
\matrix
a & a & 0 \\ 
a & a & 0 \\ 
\beta -1 & \alpha +1 & 0
\endmatrix
\right) ,\quad \nabla _{X_{3}}=\left( 
\matrix
0 & 0 & 0 \\ 
0 & 0 & 0 \\ 
0 & 0 & 0
\endmatrix
\right) 
$$
The affine representation is written 
$$
\left( 
\matrix
a(x_{1}+x_{2}) & a(x_{1}+x_{2}) & 0 & x_{1} \\ 
a(x_{1}+x_{2}) & a(x_{1}+x_{2}) & 0 & x_{2} \\ 
\alpha x_{1}+(\beta -1)x_{2} & \beta x_{1}+(\alpha +1)x_{2} & 0 & x_{3} \\ 
0 & 0 & 0 & 0
\endmatrix
\right) 
$$
\endproclaim

BIBLIOGRAPHY

\noindent

\ref \no1 \by Besnoit Y. \paper Une nilvari\'{e}t\'{e} non affine. \jour J.Diff.Geom. \vol 41 \pages 21--52
\yr 1995 \endref

\ref \no2 \by Burde D. \paper Left invariant affine structure on reductive Lie
groups. \jour J. Algebra \vol 181 \pages 884--902 \yr 1996 \endref

\ref \no3 \by Burde D. \paper Affine structures on nilmanifolds. \jour Int. J.
of Math, \vol 7 \pages 599--616 \yr 1996 \endref

\ref \no4 \by Burde D.\paper Simple left-symmetric algebras with solvable Lie
algebra. \jour Manuscripta math.\vol 95 \pages 397--411 \yr 1998 \endref

\ref \no5 \by Helmstetter J. \paper Radical d'une alg\`{e}bre sym\'{e}trique \`{a}
gauche. \jour Ann. Inst. Fourier \vol 29 \pages 17--35 \yr 1979 \endref

\ref \no6 \by Goze M., Khakimdjanov Y.
\book Nilpotent Lie Algebras \bookinfo Mathematics and its Applications 361
 \publ Kluwers Academic Publishers
 \yr 1996 \endref

\ref \no7 \by Goze M., Remm E. \paper Sur les nilvari\'{e}t\'{e}s affines. 
\inbook Actes Colloque de Brasov. \publ Univ. Brasov \yr 1999 \endref

\smallskip
Universit\'e de Haute Alsace.
Facult\'e des Sciences et Techniques.
4, rue des Fr\`eres Lumi\`ere. F.68093 MULHOUSE Cedex. France.

\smallskip 
{\it Email adress }: E.Remm\@univ-mulhouse.fr ;  M.Goze\@univ-mulhouse.fr

\end{document}